\input eplain
\input amstex

\input amssym\input wasyfont\documentstyle{amsppt}

\def\e{\epsilon}
\def\a{\alpha}
\def\b{\beta}

\def\g{\gamma}\def\k{\kappa}
\def\d{\delta}
\def\D{\Delta}
\def\s{\sigma}

\def\l{\lambda}
\def\x{\times}

\def\C{\Cal C}
\def\o{\overline}
\def\f{\flushpar}
\def\u{\underline}
\def\v{\varphi}

\def\om{\omega}
\def\Om{\Omega}
\def\B{\Cal B}\def\A{\Cal A}
\def\T{\widehat T}
\def\({\biggl(}
\def\){\biggr)}
\def\eppt{\text{ergodic, probability preserving transformation}}\def\enst{\text{ergodic, nonsingular transformation}}
\def\ppt{\text{probability preserving transformation}}

\def\wrt{\text{with respect to}}\def\st{\text{such that}}
\def\<{\langle}
\def\>{\rangle}

\def\bul{\smallskip\f$\bullet\ \ \ $}

\def\sms{\smallskip\f}\def\Par{\smallskip\f\P}
\def\pf{\smallskip\f{\it Proof}\ \ \ \ }\def\sms{\smallskip\f}\def\sbul{\f$\bullet\
\ \ $}\def\sms{\smallskip\f}
\def\lra{\longrightarrow}\def\lla{\longleftarrow}
\document\topmatter \title{Predictability, entropy and information of infinite transformations   }\endtitle

\author Jon. Aaronson $\&$ Kyewon Koh Park  \endauthor\address[Aaronson]{\ \ School of Math. Sciences, Tel Aviv University,
69978 Tel Aviv, Israel.}
\endaddress
\email{aaro\@ tau.ac.il}\endemail\address[Park]{\ \
Dept. of Math., Ajou University, Suwon 442-729,
South Korea.}
\endaddress
\email{kkpark\@ ajou.ac.kr}\endemail\abstract We  show that  a
certain type of  quasi finite, conservative, ergodic ,  measure preserving transformation  always has a  maximal zero entropy
factor, generated by predictable sets. We also construct a  conservative, ergodic ,  measure preserving transformation which is not quasi finite; and consider distribution
asymptotics  of information showing that  e.g. for Boole's transformation,
information is  asymptotically mod-normal with normalization $\propto\sqrt n$. Lastly we see that
certain  ergodic, probability preserving transformations with zero entropy have analogous properties and consequently entropy dimension of at most $\frac12$.

\endabstract\thanks{\copyright 2007 Preliminary
version.}\endthanks\subjclass 37A40,
 60F05)\endsubjclass\keywords{measure preserving transformation; conservative; ergodic; entropy; quasi finite; predictable set; entropy dimension}\endkeywords
 \endtopmatter\rightheadtext{Entropy } \heading{\tt \S0 Introduction}\endheading
Let $(X,\B,m,T)$ be a  conservative, ergodic ,   measure preserving transformation  and let $\Cal F:=\{F\in\B:\ m(F)<\infty\}$.

Call a set $A\in\Cal F$  {\it $T$-predictable} if it is measurable with respect to its own past in the sense that $A\in\s(\{T^{-n}A:\ n\ge 1\})$ (the $\s$-algebra generated by $\{T^{-n}A:\ n\ge 1\}$) and let $\Cal P=\Cal P_T:=\{\text{\tt
$T$-predictable sets}\}$.

If $m(X)<\infty$,    Pinsker's theorem ([Pi]) says that  \bul $\Cal P_T$ is the {\tt maximal, zero-entropy  factor algebra}

i.e.
 $\Cal
P\subset\B$ is a factor algebra ($T$-invariant, sub-$\s$-algebra),
$h(T,\Cal P)=0$ (see \S1) and if $\Cal C\subset\B$ is a factor algebra, with
$h(T,\Cal C)=0$, then $\Cal C\subseteq \Cal P$.  $\Cal P$ is aka the
{\it Pinsker algebra} of $(X,\B,m,T)$.

When $(X,\B,m,T)$ is a  conservative, ergodic ,   measure preserving transformation
with $m(X)=\infty$,  the above statement fails and indeed $\s(\Cal P)=\B$: Krengel  has shown ([K2]) that: \bul $\forall\ A\in\Cal F,\
\e>0,\ \exists\ B\in\Cal F,\ m(A\D B)<\e$, a {\it  strong generator}
in the sense that $\s(\{T^{-n}B:\ n\ge 1\})=\B$, whence $\s(\Cal P_T)=\B$.

It is not known if there is always a maximal,
zero-entropy factor algebra (in case there is some zero-entropy
factor algebra).

We recall the basic properties of {\tt entropy} in \S1 and define the class of {\it log lower bounded}  conservative, ergodic , measure preserving transformations in
\S2.

\par These are quasi finite in the sense of [K1] and are discussed in \S2 in this context  where also examples are constructed including a  conservative, ergodic,  measure preserving transformation  which is not quasi finite. \par A  log lower bounded  conservative, ergodic ,  measure preserving transformation   with some zero-entropy factor algebra has a maximal, zero-entropy factor algebra generated by a specified
hereditary subring of predictable sets (see \S5). \par We  obtain information convergence (in \S4)
for quasi finite transformations (cf [KS]). \par For quasi finite, pointwise dual
ergodic transformations  with regularly varying return
sequences,  we obtain (in \S6) distributional convergence of information.
Lastly, we construct a probability preserving transformation with
zero entropy with analogous distributional properties and estimate its {\tt entropy dimension} in the sense of [FP]. This example is unusual in that it has a generator with information function  asymptotic to a non degenerate random variable (the {\tt range of Brownian motion}). \heading{\tt
\S1 Entropy }\endheading

We recall the basic entropy theory of a $\ppt$ $(\Om,\A,P,S)$.  Let
$\a\subset\A$ be a countable partition.\bul The {\it entropy} of
$\a$ is $H(\a):=\sum_{a\in\a}P(a)\log\tfrac1{P(a)}$;\bul the
$S$-{\it join of $\a$ from $k$ to $\ell$} (for $k<\ell$) is
$$\a_k^{\ell}(S):=\{\bigcap_{j=k}^\ell S^{-j}a_j:\
a_k,a_{k+1},\dots,a_\ell\in\a\}.$$ \bul By subadditivity, $\exists\
\lim_{n\to\infty}\tfrac1nH(\a_0^{n-1}(S))=: h(S,\a)$ (the {\it
entropy \footnote{mean entropy rate} of $S$ $\wrt$ $\a$}). \bul The
{\it entropy of $S$} $\wrt$ the {\it factor algebra} ($S$-invariant,
$\s$-algebra), $\C\subset\A$ is
$h(S,\C):=\sup_{\a\subset\C}h(S,\a).$ \bul  By the generator
theorem, if $\a$ is a partition, then $h(S,\a)=h(S,\s(\{S^n\a:\
n\in\Bbb Z\}))$.\bul The {\it information} of the countable partition $\a\subset\A$ is the function $I(\a):\Om\to\Bbb R$ defined by
$$I(\a)(x):=\log\tfrac1{P(\a(x))}$$ where $\a(x)\in\a$ is defined by $x\in\a(x)\in\a$. Evidently
$$H(\a)=\int_\Om I(\a)dP.$$\bul {\bf  Convergence  of information} is given by
the celebrated Shannon- McMillan-Breiman theorem (see [S], [M], [Br]
respectively), the statement ($\goth I$) here being due to Chung [C] (see also [IT]).

Let $(\Om,\A,P,S)$ be an  ergodic $\ppt$ and let $\a$ be a partition
with $H(\a)<\infty$, then
$$\tfrac1nI(\a_1^N(S))\ \lra\ h(S,\a)\ \ \text{a.s. as $n\to\infty$};\tag{$\goth I$}$$
equivalently
${P(\a_1^N(S)(x))}=e^{-nh(S,\a)(1+o(1))}$ for a.e. $x\in \Om$ as $n\to\infty$
  where
$x\in\a_1^N(S)(x)\in\a_1^N(s).$

 \bul We'll need {\bf Abramov's formula} for the entropy of an induced transformation of an ergodic
$\ppt$ $(\Om,\A,P,S)$:
$$h(S_A)=\tfrac1{P(A)}h(S)\ \forall\ A\in\A$$
where $S_A:A\to A$ is the {\it induced transformation}
on $A$ defined by
$$S_Ax:=S^{\v_A(x)}x,\ \v_A(x):=\min\,\{n\ge 1:\ S^nx\in A\}\ \ \ \ (x\in A).$$
\bul Abramov's formula can be  proved using convergence of
information  (see [Ab] and \S4 here).\subheading{\smc Krengel
entropy } Suppose that $(X,\B,m,T)$ is a  conservative, ergodic , measure preserving transformation  then using
Abramov's formula (as shown in [K1])
$$m(A)h(T_A)=m(B)h(T_B)\ \ \ \forall\ A,B\in\Cal F:=\{F\in\B,\ 0<m(F)<\infty\}.$$

Set $\u h(T):=m(A)h(T_A)$,\ \ \ (any $A\in\B,\ 0<m(A)<\infty$) --
the {\it Krengel\  entropy of $T$}.

More generally, the {\it Krengel entropy of $T$ $\wrt$ the factor}\
(i.e. $\s$-finite, $T$-invariant sub-$\s$algebra) $\Cal C\subset\B$
is $$\u h(T,\Cal C):=m(A)h(T_A,\Cal C\cap A)\ \ \ (A\in\C,\
0<m(A)<\infty).$$

\bul Another definition of entropy is given in [Pa].

It is shown in [Pa] that for {\it quasi finite} (see \S2 below)  conservative, ergodic ,
  measure preserving transformations, the two entropies coincide.

\heading{\tt\S2 Quasifiniteness and Log lower boundedness}\endheading

\subheading{\smc Quasifiniteness}

\sms Let $(X,\B,m,T)$ be  a  conservative, ergodic ,  measure preserving transformation .

Recall from [K1] that a set $A\in\Cal F$ is called {\it quasi
finite} (qf) if $H_A({\rho_A})<\infty$ where $\rho_A:=\{A\cap
T^{-n}A\setminus\bigcup_{j=1}^{n-1}T^{-j}A:\ n\ge 1\}$ and that $T$
is so called if $\exists$ such a set. As shown in
proposition 7.1 in [K1], \bul  for  $A\in\Cal F$ quasi finite, $A\in\Cal P_T\
\iff\ h(T_A,\rho_A)=0$.

There are  conservative, ergodic ,  measure preserving transformation s which are not quasi finite. An unpublished example
of such by Ornstein is mentioned in [K2, p. 82].

Here we construct a  conservative, ergodic ,  measure preserving transformation  with no quasi finite extension. To do
this we first establish a saturation property for the collection of
quasi finite sets: \proclaim{Proposition 2.0}

 Suppose that $(X,\B,m,T)$ is a  conservative, ergodic , quasi finite,  measure preserving transformation , then
$\forall\ F\in\Cal F,\ \exists\ A\in\B\cap F$ $\st\ m(A)>0$ and
$\st$ each $B\in\B\cap A$ is quasi finite.\endproclaim \pf We show first
that \Par1 {\sl if $F\in\Cal F$ is quasi finite, then $\forall\ \e>0,\
\exists\ A\in\B\cap F$  $\st\ m(F\setminus A)<\e$ and $\st$ each
$B\in\B\cap A$ is quasi finite. }\pf By  ($\goth I$),
$\tfrac1nI(\rho_F)_0^{n-1}(T_F))\to h(T_F,\rho_F)$
a.e. as $n\to\infty$. By Egorov's  theorem, $\exists\ A\in\B\cap F\
\st\ m(F\setminus A)<\e$ and $\st$ the convergence is uniform on
$A$.

For $B\in\B\cap A$,  let $N_{n,B}:=\#\,\{a\in(\rho_F)_0^{n-1}(T_F):\
m(a\cap B)>0\}$ (where $\# F$ means the number of elements in the set $F$), then $N_{n,B}=e^{nh(T_F,\rho_F)(1+o(1))}$ as
$n\to\infty$.

Define $\psi:B\to\Bbb N$ by $\psi(x):=\min\,\{n\ge 1:\ T_F^nx\in
B\}$, then \sms - $\int_B\psi dm=\sum_{n=1}^\infty
nm([\psi=n])=m(F)<\infty$ (by Kac's formula);\sms -
$\v_B(x)=\sum_{j=0}^{\psi(x)-1}\v_F(T_F^jx)$ whence
$$\rho_B\prec\g_B:=\bigcup_{n=1}^\infty\{[\psi=n]\cap a:\
a\in(\rho_F)_0^{n-1}(T_F)\}.$$ Thus
$$\align H_{m_B}(\rho_B) &\le H_{m_{B}}(\g_B)\\ &=\sum_{n=1}^\infty m_B([\psi=n])H_{m_{[\psi=n]}}((\rho_F)_0^{n-1}(T_F))\\ &\le
\sum_{n=1}^\infty m_B([\psi=n])\log N_{n,B}<\infty\ \because\ \log
N_{n,B}\sim nh(T_F,\rho_F).\checkmark\qed\P1\endalign$$ To complete the proof,
let $F\in\Cal F$. Suppose that $Q\in\Cal F$ is quasi finite, then evidently
so is $T^{-n}Q\ \forall\ n\ge 1$. By ergodicity, $\exists\ n\ge 1\
\st\  m(F\cap T^{-n}Q)>0$. By \P1, $\exists\ G\in\B\cap T^{-n}Q$
$\st\ m(T^{-n}Q\setminus G)<\e:=\tfrac{m(F\cap T^{-n}Q)}9$ and $\st$
each $B\in\B\cap G$ is quasi finite. The set $A=G\cap F$ is as
required.\qed \subheading{Example 2.1}\ \ \ \ \ \ Let
$(X_0,\B_0,m_0,T_0)$ be the  conservative, ergodic ,  measure preserving transformation  defined as in [Fr] by
the cutting and stacking construction
$$B_0=1,\ B_{n}=\bigoplus_{k=1}^{N_{n}}B_{n-1}0^{L_{n,k}}$$
where $N_n,\ L_{n,k}\ 1\le k\le N_n$  satisfy
$$N_{n+1}\ge e^{nN_1\dots N_n},\ \ L_{n,k+1}> \sum_{j=1}^kL_{n,j}+kh_{n-1},$$ where $h_n:=|B_n|$.
\proclaim{Proposition 2.1}\sms No extension $T$ of the  conservative, ergodic ,  measure preserving transformation
$T_0$ defined in example 2.1 is  quasi finite.\endproclaim\pf Suppose
otherwise, that $(X,\B,m,T)$ is a (WLOG)  conservative, ergodic  extension of $T_0$
and that $F\in\Cal F$ is quasi finite, then evidently so is $T^nF\ \forall\
n\ge 1$. By  proposition 2.0 $\exists\ A\in\B,\
A\subset B_0$ quasi finite. We'll contradict this (and therefore the
assumption that $\exists\ F\in\Cal F$  quasi finite). \Par1 Write
$B_n=\bigcup_{j=0}^{h_n-1}T^jb_n$ where $b_n\subset B_0,\
m(b_n)=\tfrac1{N_1N_2\dots N_n}$ and
$B_n=\biguplus_{k=1}^{N_{n+1}}B_n^{(k)}=\biguplus_{k=1}^{N_{n+1}}T^{\k(n+1,k)}B_n^{(1)}$
where $\k(n+1,1)=0$ and
$\k(n+1,k)=(k-1)|B_n|+\sum_{j=1}^{k-1}L_{n+1,j}$ (i.e. the
$B_n^{(k)}\ \ (1\le k\le N_{n+1})$ are the subcolumns of $B_n$
appearing in $B_{n+1}$). \Par2 For $n\ge 1$, let $\goth k_n:=\{0\le
j\le h_n-1:\ T^jb_n\subset B_0\}$, then
 \sms  $B_0=\biguplus_{j\in\goth k_n}T^jb_n$, $|\goth k_n|=N_1N_2\dots N_{n}$ and \sms  for $x\in b_{n}$,
$\{T_{B_0}^kx\}_{k=0}^{N_1N_2\dots N_{n}-1} = \{T^jx:\ j\in \goth
k_n\}$.

\Par3 Fix $0<\e<\tfrac13$ and let
$$b_{n,\e}:=\{x\in b_{n+1}:\ |\tfrac1{|\goth k_{n+1}|}\sum_{k\in\goth k_{n+1}}1_A(T^kx)-m(A)|< \e m(A)\}.$$
By \P2 above,  for $x\in b_{n+1}$,
$$\tfrac1{|\goth k_{n+1}|}\sum_{k\in\goth k_{n+1}}1_A(T^kx)=\tfrac1{N_1N_2\dots N_{n+1}}\sum_{k=0}^{N_1N_2\dots N_{n}-1}1_A(T_{B_0}^kx)$$ and a standard argument using the ergodic theorem for $T_{B_0}$ shows that  $\exists\ M$ so that
$m(b_{n,\e})>(1-\e)m(b_{n+1})\ \ \forall\ n\ge M.$
 \Par4
Fix   $n\ge M$ and $x\in b_{n+1}$, let $\goth k_{A,n,x}:=\{k\in\goth
k_{n+1}:\ T^kx\in A\}$ and $A_{n,x}:=\{T^jx\}_{j\in \goth
k_{A,n,x}}$, then for $x\in b_{n,\e}$,
$$\#\{1\le k\le N_{n+1}:\
A_{n,x}\cap B_n^{(k)}\ne\phi\}\ge (1-e)m(A)\tfrac{|\goth
k_{n+1}|}{h_n}=(1-\e)m(A)N_{n+1}.$$ \bul For $n\ge M,\ x\in
b_{n,\e}$,  write
$$\{1\le k\le N_{n+1}:\
A_{n,x}\cap B_n^{(k)}\ne\phi\}=:\{\k_i(x):\ 1\le i\le \nu\}$$ where
$\nu-1>(1-\e)m(A)N_{n+1}$ and $\k_i(x)<\k_{i+1}(x)\ \forall\ i$.
\bul For $1\le i\le \nu$, let
 $\goth k_{A,n,x}^{(i)}:=\{k\in\goth k_{n+1}:\ T^kb_{n+1}\subset A_{n,x}\cap B_n^{(\k_i)}\}$ and let $\u m_i:=\min\,\goth k_{A,n,x}^{(i)},\ \o m_i:=\max\,\goth k_{A,n,x}^{(i)}$;
 $y_i:=\o m_{i+1}-\u m_i,\ \ (1\le i\le \nu-1)$. Note that
 \sms  $y_i\le \sum_{j=1}^{\k_i}L_{n+1,j}+\k_ih_n<L(n+1,\k_i+1)\le L(n+1,\k_{i+1})\le y_{i+1}$.
\Par5 For $K\subset\goth k_{n+1}$, let $a_K:=\{x\in b_{n+1},\ \goth
k_{A,n,x}=K\}$ and let $$\b_n:=\{a_K:\ K\subset\goth k_{n+1}\},\ \
\a_n:=\{\widehat a:=\bigcup_{j\in \goth k_n}T^ja:\ a\in\b_n\}.$$
 \bul For  $a\in\b_n,\ a\subset b_{n,\e},\ 1\le i\le \nu-1$,\sms
$m(a\cap [\v_A=y_i(a)])=\tfrac{m(a)}{N_1\dots N_{n+1}}.$ \bul Thus
$$\align H(\rho_A) & \ge H(\rho_A\|\a_n)\\ &\ge\sum_{a\in\b_n,\ a\subset
b_{n,\e}}m(a)
\sum_{i=1}^{\nu-1}m([\v_A=y_i(a)]|a)\log\tfrac1{m([\v_A=y_i(a)]|a)}\\
&\ge m(\widehat b_{n,\e})\tfrac{(\nu-1)\log (N_{n+1})}{N_1\dots N_{n+1}}\\
&\ge (1-\e)^2m(A)\tfrac{\log N_{n+1}}{N_1N_2\dots N_{n}}\\
&>(1-\e)^2m(A)n \uparrow\infty.\qed\endalign$$

\subheading{\smc Log lower boundedness}
\sms For $(X,\B,m,T)$   a  conservative, ergodic ,  measure preserving transformation ; set
$$\Cal F_{\log,T}:=\{A\in\B:\ 0<m(A)<\infty,\ \int_A\log\v_Adm<\infty\}.$$
\bul Note that \sms  $\Cal F_{\log,T}\subset\{\text{\tt quasi finite
sets}\}$; because
$$p_n\ge 0,\ \sum_{n=1}p_n\log n<\infty\
\implies\ \sum_{n=1}p_n\log\frac1{p_n}<\infty.\tag{$\bigstar$}$$ \bul Call $T$ {\it
$\log$-lower bounded }(LLB) if $\Cal F_{\log,T}\ne\emptyset$.
\proclaim{Proposition 2.2}\sms {\rm (i)} $T$ is LLB iff
$\tfrac1{\log n}\sum_{k=0}^{n-1}f\circ T^n\to\infty$ a.e. as
$n\to\infty$ for some and hence all $f\in L^1(m)_+:=\{f\in L^1,\
f\ge 0,\ \int_Xfdm>0\};$\sms (ii) $T$ is not LLB iff
$\liminf_{n\to\infty}\tfrac1{\log n}\sum_{k=0}^{n-1}f\circ T^n=0$
a.e.  for some and hence all $f\in L^1_+;$ \sms {\rm (iii)}  If
$(X,\B,m,T)$ is LLB and $\Cal C\subset\B$ is a factor, then $\Cal
C\cap\Cal F_{\log,T}\ne\emptyset$.\sms {\rm (iv)} $\Cal F_{\log,T}$
is a hereditary ring.
\endproclaim\pf Statements (i) and (ii) follow from theorem 2.4.1 in [A] and
(iii) follows from these. We prove (iv). \sms\par Suppose that
$A\in\Cal F_{\log,T},\ B\in\B,\ B\subset A$, then
$\v_B(x)=\sum_{k=0}^{\psi(x)-1}\v_A(T_A^kx)\ \ \ (x\in B)$ where
$\psi:B\to\Bbb N,\ \psi(x):=\min\,\{n\ge 1:\ T_A^nx\in B\}$.

By Kac formula,
$$\int_B\sum_{k=0}^{\psi-1}f\circ T_A^kdm=\int_Afdm\ \forall\ f\in L^1(m).$$ To see that $B\in\Cal F_{\log,T}$,
we use this and $\log(k+\ell)\le \log(k)+\log(\ell)$:
$$\align \int_B\log\v_Bdm &=\int_B\log(\sum_{k=0}^{\psi-1}\v_A\circ T_A^k)dm\\ &\le \int_B\sum_{k=0}^{\psi-1}\log(\v_A\circ T_A^k)dm\\ &=\int_A\log\v_Adm<\infty.\endalign$$

 Suppose that $A, B\in\Cal F_{\log,T}$, then $\v_{A\cup B}\le
1_A\v_A+1_B\v_B$ whence
$$\align \int_{A\cup B}\log(\v_{A\cup B})dm &= \int_{A}\log(\v_{A\cup B})dm+\int_{B}\log(\v_{A\cup B})dm\\ &\le
\int_{A}\log(\v_A)dm+\int_{B}\log(\v_B)dm\\ &<\infty.\qed\endalign$$

\heading{\tt \S3 Examples of LLB
transformations}\endheading\subheading{\smc Pointwise dual ergodic
transformations}

A  conservative, ergodic ,  measure preserving transformation  $(X,\B,m,T)$ is called {\it pointwise dual
ergodic} if there is a sequence of constants $(a_n(T))_{n\ge 1}$
\f(called the {\it return sequence} of $T$) so that
$$\tfrac1{a_n(T)}\sum_{k=0}^{n-1}\T^k f\to\int_Xfdm\ \ \ \text{a.e. for some (and hence all)}\ \ f\in
L^1(m)_+$$ where $\T:L^1(m)\to L^1(m)$ is the {\it transfer
operator} defined by $$\int_A\T fdm=\int_{T^{-1}A}fdm\ \ \ (f\in
L^1(m),\ A\in\B).$$ See [A, 3.8].\proclaim{Proposition 3.1}\sms Let
$(X,\B,m,T)$ be a pointwise dual ergodic,  conservative, ergodic , measure preserving transformation , then $T$
is LLB $\iff$ $\sum_{n=1}^\infty\tfrac1{na_n(T)}<\infty.$
\endproclaim\pf\sms  Let
$A\in\Cal F$ be a {\it uniform set} in the sense that for some $f\in
L^1(m)_+$
$$\tfrac1{a_n(T)}\sum_{k=0}^{n-1}\T^k f\to\int_Xfdm\ \ \text{\rm uniformly on}\ \ A.$$
By lemma 3.8.5 in [A], $$\int_A(\v_A\wedge
n)dm=m(\bigcup_{k=0}^nT^{-k}A)\asymp \tfrac{n}{a_n(T)}$$ whence
$$A\in\Cal F_{\log}\ \iff\
\sum_{n=1}^\infty\tfrac{m(\bigcup_{k=0}^nT^{-k}A)}n<\infty\ \iff\
\sum_{n=1}^\infty\tfrac1{na_n(T)}<\infty.\qed$$
 \subheading{Remarks}\par 1) For example, the  simple random walk on $\Bbb Z$ is LLB
($\because\ a_n(T)\propto\sqrt n$); whereas the  simple random walk
on $\Bbb Z^2$ is not LLB ($\because\ a_n(T)\propto\log n$).\par 2)
It is not known whether the simple random walk on $\Bbb Z^2$ is
quasi finite, or even has a factor with finite entropy. \proclaim{Example
3.2}\sms There is  a quasi finite,  conservative, ergodic, Markov shift $(X,\B,m,T)$ with
\f $a_n(T)\asymp \sqrt{\log n}$.\endproclaim \bul Note that by proposition 3.1,
this $T$ is not   LLB. \pf {\it of example 3.2 }: Let $f_{4^{4^n}}:=\tfrac1{{2^n}}\ \ n\ge 1$
and $f_k:=0\ \forall\ k\in \Bbb N\setminus 4^{4^{\Bbb N}}$, then
$f\in\Cal P(\Bbb N)$.\bul Let $\Om:=\Bbb N^\Bbb Z$ and let $P=f^\Bbb
Z\in\Cal P(\Om,\B(\Om))$ be product measure, then
\f$(\Om,\B(\Om),P,S)$ in an  ergodic, probability preserving transformation  where $S:\Om\to\Om$ is the
shift. \bul Define $\v:\Om\to\Bbb N$ by $\v(\om):=\om_0$ and let
$(X,\B,m,T)$ be the tower over $(\Om,\B(\Om),P,S)$ with height
function $\v$.\bul It follows that $(X,\B,m,T)$ is a  conservative, ergodic , Markov
shift with $a_n(T)\asymp\sum_{k=0}^nu_k$ where $u$ is defined by the
renewal equation: $u_0=1,\ u_n=\sum_{k=1}^nf_ku_{n-k}$.\bul To see
that $(X,\B,m,T)$ is quasi finite, we check that $\Om$ is quasi finite. Indeed
$$H_\Om(\rho_\Om)=\sum_{k\ge 1,\ f_k>0}^\infty
f_k\log\tfrac1{f_k}=\sum_{n=1}^\infty\tfrac{n\log 2}{2^n}<\infty.$$
\bul To estimate $a_n(T)$, recall that by lemma 3.8.5 in [A],
$a_n(T)\asymp\tfrac{n}{L(n)}$ where
$$L(n):=m(\bigcup_{k=0}^nT^{-k}\Om)=\sum_{k=0}^n\sum_{\ell=k+1}^\infty f_\ell.$$
Now, $$\sum_{\ell=k+1}^\infty f_\ell=\sum_{n>\log_4\log_4
k}\tfrac1{2^n}\asymp \tfrac1{2^{\log_4\log_4
k}}=\tfrac1{\sqrt{\log_4 k}}.$$ Thus $L(n)\asymp\tfrac{n}{\sqrt{\log
n}}$ and $a_n(T)\asymp \sqrt{\log n}$.\qed \subheading{\smc The
Hajian-Ito-Kakutani transformations}\bul Let $\Om=\{0,1\}^{\Bbb N},\
\ell(\om):=\min\,\{n\ge 1:\ \om_n=0\}$ and let $\tau:\Om\to\Om$ be
the {\it adding machine} defined by
$$\tau(1,\dots,1,0,\om_{\ell(\om)+1},\dots):=(0,\dots,0,1,\om_{\ell(\om)+1},\dots).$$
For $p\in (0,1)$, define $\mu_p\in\Cal P(\Om)$ by
$\mu_p([a_1,\dots,a_n]):=p_{a_1}\dots p_{a_n}$ where $p_0:=1-p,\
p_1:=p$. It follows that $(\Om,\A,\mu_p,\tau)$ is an $\enst$ with
$\tfrac{d\mu_p\circ\tau}{d\,\mu_p}=(\tfrac{1-p}p)^\phi$ where
$\phi:=\ell-2$.

Now let $X:=\Om\x\Bbb Z$ and define $T:X\to X$ by $T(x,n)=(\tau
x,n+\phi(x))$. For $p\in (0,1)$, define $m_p\in\goth M(X)$ by
$m_p(A\x\{n\}):=\mu_p(A)(\tfrac{1-p}p)^{-n}$.

As shown in [HIK] (see also [A]) $T_p=(X,\B,m_p,T)$ is a  conservative, ergodic , measure preserving transformation
(aka  the {\it Hajian-Ito-Kakutani} transformation). The entropy is given by
$\u h(T_p)=h((T_p)_{\Om\x\{0\}})=0$ by [MP] since $(T_p)_{\Om\x\{0\}}$ is
the Pascal adic transformation.\proclaim{Proposition 3.3}\par $(X,\B,m_p,T)$ is LLB $\forall \ 0<p<1$.
\endproclaim\pf  As in the proof of proposition 5.1 in [A1],
$$\align \sum_{k=0}^{2^n-1}1_{\Om\x\{0\}}\circ T^k(x,0) &=
\#\{0\le k\le 2^n-1:\ \sum_{j=0}^{k-1}\phi(\tau^jx)=0\}\\ &\ge
\#\{0\le K\le n-1:\ \sum_{j=0}^{2^K-1}\phi(\tau^jx)=0\}\endalign$$
Now $\sum_{j=0}^{2^K-1}\phi(\tau^jx)=\phi(S^Kx)$ where $S:\Om\to\Om$ is the shift, and so
$$\sum_{k=0}^{2^n-1}1_{\Om\x\{0\}}\circ T^k(x,0)\ge \#\{0\le K\le n-1:\ \phi(S^Kx)=0\}
\sim (1-p)n$$
 for $\mu_p$-a.e. $x\in\Om$ by Birkhoff's theorem for the ergodic, probability preserving transformation  $(\Om,\B(\Om),\mu_p,S)$.
 The LLB property now follows from proposition 2.2.\qed

\bul Let $\goth G$ be the Polish group of  measure preserving transformation s of \f$(\Bbb
R,\B(\Bbb R),m_{\Bbb R})$ equipped with the weak topology.
\proclaim{Proposition 3.4} \sms The collection of LLB  measure preserving transformation s is
meagre in $\goth G$. \endproclaim\pf Let
 $$\text{\it\$}:=\{T\in\goth G:\ \exists\ n_k\to\infty,\
\tfrac{S_{n_k}(f)}{\log n_k}\to 0\ \text{\rm a.e.\ }\ \forall\ f\in
L^1\}$$ where $S_n(f)=S_n^T(f):=\sum_{j=0}^{n-1}f\circ T^j$. \sms By
proposition 2.2, it suffices to show that {\it\$} is a dense $G_\d$
set in $\goth G$.
\par By example 3.2, $\exists$ a  conservative, ergodic ,  measure preserving transformation  $T\in\text{\it\$}$.
$\text{\it\$}$ is conjugacy invariant, and so dense in $\goth G$ by
the conjugacy lemma (e.g. 3.5.2 in [A]).
\smallskip\sms
\par To see that {\it\$} is a $G_\d$ set, let \bul $P\sim m$ be a
probability;\bul fix $\{A_n:\ n\in\Bbb N\}\subset\Cal F:=\{A\in\B:m(A)<\infty\}$ so that $\s(\{A_n:\ n\in\Bbb N\})=\B$ and  let
$$\text{\it\$}':=\bigcap_{k=1}^\infty\bigcup_{n=k}^\infty\bigcap_{\nu=1}^k
\{T\in\goth G:\ P([S_n(1_{A_\nu})>\tfrac1k\log n]) <{1\over 2^k}
\},$$ then $\text{\it\$}'$ is a $G_\d$. We claim
$\text{\it\$}'=\text{\it\$}$.
\smallskip
Evidently,
$$\text{\it\$}'=\{T\in\goth G:\ \exists\ n_k\to\infty\text{ such that }
\tfrac{S_{n_k}(1_{A_\nu})}{\log n_k}\to 0 \text{  a.e. }\forall\
\nu\ge 1\}$$ whence $\text{\it\$}'\supset\text{\it\$}$.
\par Now suppose that $T\in\text{\it\$}'$, that $\tfrac{S_{n_k}(1_{A_\nu})}{\log n_k}\to 0$
  a.e. $\forall\ \nu\ge 1$ and let $f\in L^1$. Evidently
$\tfrac{S_{n}(f)}{\log n}\to 0$
  a.e. on $\goth D$, the dissipative part of $T$. The conservative part of $T$ is
  $$\goth C=\bigcup_{\nu=1}^\infty\hat A_\nu\ \text{\tt where}\ \hat A_\nu:=[
\sum_{n=1}^\infty 1_{A_\nu}\circ T^n=\infty].$$
 By Hopf's theorem,
${S_n(f)(x)\over S_n(1_{A_\nu})(x)}\to h_\nu(f)$ a.e. on $A_\nu\
\forall\ \nu\ge 1$ where $h_\nu(f)\circ T=h_\nu(f)$ and
$\int_{A_\nu}h_\nu(f)dm=\int_Xfdm$, whence, a.e. on $\hat A_\nu$,
$$\tfrac{S_{n_k}(f)}{\log n_k}=\tfrac{S_{n_k}(f)}{S_{n_k}(1_{A_\nu})}\cdot\tfrac{S_{n_k}(1_{A_\nu})}{\log n_k}\to 0.\qed$$  \heading{\tt\S4  Information convergence
}\endheading \sms Let $(X,\B,m,T)$ be a  conservative, ergodic ,  measure preserving transformation . \bul  A
countable partition $\xi\subset\B$  is called {\it cofinite} if
$\exists\ A=A_\xi\in\Cal F$ with $A^c\in\xi$. We call  $A^c$  the
{\it cofinite atom} of $\xi$ and $A$ the {\it (finite) core} of
$\xi$.. \bul If $\xi\subset\B$ is cofinite, then $\xi_k^\ell(T)$ is
also cofinite, with core $A_{\xi_k^\ell(T)}=\bigcup_{j=k}^\ell
T^{-j}A$. \sms\par The $T$-process generated by a cofinite partition
$\xi$ restricted to its core $A$ is given by \sms {\bf Krengel's
formula} [K1]:
$$\xi_1^{\v_n(x)} (T)(x)\ =\ (\rho_A\vee((\xi\cap
A)\vee\rho_A)_1^n(T_A))(x)\text{\rm\  for a.e. $x\in A$}\tag{$\goth K$}$$ where
 for $x\in X,\ \a$ a partition of $X$, $\a(x)$ is defined by
$x\in\a(x)\in\a$; \sms $\v_n(x):=\sum_{k=0}^{n-1}\v_A(T_A^kx)$; and
 \sms $\rho_A:=\{A\cap T^{-n}A\setminus\bigcup_{k=1}^{n-1}T^{-k}A:\ n\in\Bbb N\}$.\bul A
cofinite partition $\xi\subset\B$ is called {\it quasi-finite} (qf)
if  $A=A_\xi$ is quasi finite and $H_A(\xi)<\infty$. \bul Note that $\xi$
quasi finite $\Rightarrow$ $H_{A}(\xi\vee\rho_{A}\vee
T_{A}\,\rho_{A})<\infty$.

\subheading{\smc Convergence of information for quasi finite partitions}

\proclaim{Proposition 4.1\ \ (c.f.
[KS])}

Let $(X,\B,m,T)$ be a  conservative, ergodic ,  measure preserving transformation , let $\xi\subset\B$ be a quasi finite
partition and let $p\in L^1(m),\ p>0,\ \int_Xpdm=1$, then for a.e.
$x\in X$,
$$\tfrac1{S_n(p)(x)}I(\xi_1^n(T))(x)\to \u h(T,\xi)$$
where $S_n(p)(x):=\sum_{k=0}^{n-1}p(T^kx)$ and $I(\xi_1^n(T))(x):=\log\tfrac1{m(\xi_1^n(T)(x))}$.\endproclaim\pf Let $A$
be the core of $\xi$ and set  $\varsigma:=(\xi\cap A)\vee\rho_A$,
then { by ($\goth K$)}
$$\varsigma_0^{s_n(x)}(T_A)(x)\subseteq\xi_1^n(T)(x)\subseteq\varsigma_1^{s_n(x)-1}(T_A)(x)\ \
\text{\rm a.e. $x\in A$}$$ where $x\in\xi(x)\in\xi$,
$s_n:=S_n(1_A)$. \sbul By ($\goth I$), for $T_A$,  a.e. on $A$,
$I(\varsigma_1^{N}(T_A))\sim Nh(T_A,\varsigma)$,\
whence for a.e.  $x\in A$,
$$\align \log\tfrac1{m(\xi_1^n(T)(x))} &\sim
\log\tfrac1{m(\varsigma_1^{s_n(x)}(T_A)(x))}\\ &\sim
s_n(x)h(T_A,\varsigma)\\ &\sim S_n(p)(x)m(A)h(T_A,\varsigma)\\
&=S_n(p)(x)\u h(T,\xi).\endalign$$ We obtain convergence  a.e. on
$\bigcup_{k=0}^NT^{-k}A$ by substituting $\xi_1^N(T)$ for $\xi$;
whence convergence a.e. on $X$
 as $\bigcup_{k=0}^NT^{-k}A\uparrow X$.\qed
\bul  Abramov's formula is proved analogously in case $(X,\B,m,T)$
is an  ergodic, probability preserving transformation. As in [Ab]:
$$h(T,\xi)\overset{\text{\rm ($\goth I$)}}\to\lla\tfrac1n\log\tfrac1{m(\xi_1^n(T)(x))}\approx
\tfrac1n s_n(x)h(T_A,\varsigma)\overset{\text{\rm Birkhoff's
PET}}\to{\lra} m(A)h(T_A,\varsigma).$$

\heading{\tt\S5 Pinsker algebra}\endheading \sms Let
$(X,\B,m,T)$ be a  LLB,  conservative, ergodic , measure preserving transformation .

Define $$\Cal F_{\Pi}:=\{A\in\Cal F_{\log,T}:\ A\in\s(\{T^{-k}A:\
k\ge 1\})\}=\Cal P\cap\Cal F_{\log,T}.$$ In this section, we show that (in case $\Cal F_{\Pi}\ne\emptyset$) $\B_{\Pi}:=\s(\Cal F_{\Pi})$ is the maximal zero entropy factor of $T$.

To do this, we'll need \proclaim{Krengel's predictability lemma}\
[K1]:\ \sms Let $(X,\B,m,T)$ be a quasi finite,   conservative, ergodic , measure preserving transformation , let
$\xi\subset\B$ be a quasi finite partition with core $A$, and let
$\zeta=\xi\cap A$, then
$$\xi\subset\xi_1^\infty(T)\ \mod m\ \iff\ h(T_A,\zeta\vee\rho_A)=0.$$ In particular
$$A\in\s(\{T^{-n}A:\ n\ge 1\})\ \iff \
h(T_A,\rho_A)=0.$$\endproclaim \bul For $F\in\Cal
F$, set
$$\Cal P_F=\Cal P_{T_F}:=\{A\in\B\cap F:\ A\in\s(\{T_{F}^{-k}A:\ k\ge 1\})\}.$$

  By Pinsker's theorem ([Pi]), \sbul $\Cal P_F$ is a $T_{F}$-factor algebra
of subsets of $F$, $h(T_{F},\Cal P_F)=0$ and \sbul if $\A\subset\B\cap
F$ is another $T_{F}$-factor algebra of subsets of $F$ with
$h(T_{F},\A)=0$, then $\A\subset \Cal P_F$. \proclaim{Theorem
5.1} \sms {\rm (i)} $\Cal F_{\Pi}$ is a ring and $\Cal F_{\Pi}\cap F=\Cal P_F\
\forall\ F\in\Cal F_{\Pi}$.

\sms {\rm (ii)} If $\Cal F_{\Pi}\ne\emptyset$, then $\s(\Cal F_{\Pi})$ is the maximal factor of zero entropy.
\endproclaim\subheading{\tt Proof}

\Par1 Let $A\in\Cal F_{\log}$. By Krengel's predictability lemma,   $F\in\Cal
F_{\Pi}$ iff $h(T_F,\rho_F)=0$.

 Thus, $F\in\Cal F_{\Pi}$ iff $\exists$ a factor $\B_0$ with
$F\in\B_0$ and $\u h(T,\B_0)=0$.

\Par2 Next, fix $F\in\Cal F_{\Pi}$. We claim  that
$\rho_F\subseteq\Cal P_F$. This is because $F\in\Cal F_{\Pi}\
\Rightarrow\ h(T_F,\rho_F)=0$. \sms\P3 We now show that $\Cal
P_F\subseteq\Cal F_{\Pi}\cap F\ \forall\ F\in\Cal F_{\Pi}$.
\sms{\it Proof}: Fix $F\in\Cal F_{\Pi}$ and let $\B_0:=\s\{T^nA:\
n\in\Bbb Z,\ A\in\Cal P_F\}$, then $\B_0$ is a factor, $F\in\B_0$
and $\B_0\cap F=\Cal P_F$. Thus $\u h(T,\B_0)=h(T_F,\Cal P_F)=0$
and by \P1 $\Cal P_F\subseteq\Cal F_{\Pi}\cap F$.\smiley
 \sms\P4  Now we claim that  $A,B\in\Cal F_{\Pi}\
\Rightarrow\ A\cup B\in\Cal P_{F}$.

\sms{\it Proof}: Set $C:=A\cup B$, then $C{\in}\Cal F_{\log,T}$. Set
$\zeta:=\{A\cap B,A\setminus B, B\setminus A\}$ and
$\xi:=\zeta\cup\{C^c\}.$ By ($\goth K$), $$\xi_1^\infty(T)\cap
C=\rho_C\vee(\zeta\vee\rho_C)_1^\infty(T_C).$$ By assumption,
$\zeta\subset\xi_1^\infty(T)\cap C$, whence also $\rho_C\subset
\xi_1^\infty(T)\cap C$. Thus $$\zeta\vee\rho_C\subset
\rho_C\vee(\zeta\vee\rho_C)_1^\infty(T_C);\ \ \ \ \therefore\ \
\zeta\vee\rho_C\vee T_C\,\rho_C\subset
 (\zeta\vee\rho_C\vee T\rho_C)_1^\infty(T_C),$$ and (using
$H_C(\zeta\vee\rho_C\vee T_C\,\rho_C)<\infty$) we have
$$h(T_C,\rho_C)\le h(T_C,\zeta\vee\rho_C\vee T_C\,\rho_C)=0$$ whence
$C\in\s(\{T^{-k}C:\ k\ge 1\})$ and $C\in\Cal F_{\Pi}$. \smiley
\sms\P5 Now we show that $\Cal F_{\Pi}$ is a ring by proving that
$A,B\in\Cal F_{\Pi}\ \Rightarrow\ \zeta:=\{A\cap B,A\setminus B,
B\setminus A\}\subset\Cal F_{\Pi}$.  \demo{Proof} By \P3, it
suffices to show that $\zeta\subset\Cal P_{C}$ where
 $C:=A\cup B$.  To see this, fix $a\in\zeta$, then
 $$h(T_C,\{a,C\setminus a\})\le h(T_C,\zeta)\le h(T_C,\zeta\vee\rho_C\vee T_C\,\rho_C) =0$$ (as above) and $a\in\Cal P_{C}$.
\smiley\enddemo  \sms\P6 To complete the proof of (i), we show that $\Cal F_{\Pi}\cap
F\subseteq\Cal P_F\ \forall\ F\in\Cal F_{\Pi}$.  \sms{\it Proof}:
Fix $F\in\Cal F_{\Pi},\ A\in\Cal F_{\Pi}\cap F$. Let
$\zeta:=\{A,F\setminus A\},\ \xi:=\zeta\cup\{F^c\}$.

By the ring property, $A\in\Cal F_{\Pi}$, whence
$\xi\subset\xi_1^\infty(T)\ \mod m$. By proposition 4,
$h(T_F,\zeta\vee\rho_F)=0$, whence
$$h(T_F,\zeta)\le h(T_F,\zeta\vee\rho_F)=0$$ and $A\in\Cal
P_F$.\smiley\Par7 To see (ii), fix $F\in\Cal F_{\Pi}$, then by (i),
$\Cal F_{\Pi}\cap F=\Cal P_F=\Cal F_{\Pi}\cap F\cap F$ whence
$\u h(T,\s(\Cal F_{\Pi}))=m(F)h(T_F,\Cal P_F)=0$ and if $\C\subset\B$ is a factor with $\u h(T,\C)=0$, then by \P1,
$\C\cap\Cal F_{\log}\subset\Cal F_\Pi$, whence $\C\subset\s(\Cal F_\Pi)$.\qed

\heading{\tt\S6 Asymptotic distribution of information with infinite invariant measure}\endheading
\subheading{\smc Pointwise dual ergodic transformations} \sms Let
$(X,\B,m,T)$ be a pointwise dual ergodic  measure preserving transformation  and assume that the
return sequence $a_n=a_n(T)$ is regularly varying with index $\a\ \
(\a\in [0,1])$, then by the Darling-Kac theorem (theorem 3.6.4 in
[A] -- see also references therein),
$$\tfrac1{a_n}S_n^T(f)\overset{\goth
d}\to\rightarrow \int_Xfdm\cdot X_\a\ \ \text{\tt as}\ n\to\infty\ \forall\ f\in L^1(m)_+\tag{\mercury}$$ where \bul $X_\a$ is a Mittag-Leffler random variable of order $\a$ normalised so that $E(X_\a)=1$; and
\bul $F_n\overset{\goth
d}\to\rightarrow Y$ means
$$\int_XG(F_n)dP \to
E(G(Y))\ \forall\ P\in\Cal P(X,\B),\ P\ll m,\ \ G\in
C([0,\infty]).$$ \sms Note that $X_1\equiv 1$,  $X_0$
has exponential distribution and for $\a\in
(0,1),\ X_\a=\tfrac1{Y_\a^\a}$ where $E(e^{-tY_\a})=e^{-ct^\a}$
(some $c=c_\a>0$). In particular $X_{\frac12}=|\Cal N|$ where $\Cal
N$ is a centered Gaussian random variable on $\Bbb R$. \proclaim{Proposition 6.1} Suppose that $(X,\B,m,T)$ is a
quasi finite, pointwise dual ergodic  measure preserving transformation   and assume that the return
sequence $a_n=a_n(T)$ is regularly varying with index $\a\ \ (\a\in
[0,1])$. \par If $\xi\subset\B$ is quasi finite, then
$$\tfrac1{a_n(T)}\log\tfrac1{m(\xi_1^n(T)(x))}\overset{\goth
d}\to\rightarrow \u h(T,\xi)X_\a$$ as $n\to\infty$.\endproclaim\pf
This follows from proposition 4.1 and (\mercury).\qed

 \subheading{Example 6.2:\ \ \ {\tt Boole's transformation}}
\sms Let $(X,\B,m,T)$ be given by $X=\Bbb R,\ m=$ Lebesgue measure
and $Tx=x-\tfrac1x$, then $T$ ( see [A]) is a pointwise dual
ergodic,  measure preserving transformation  with $a_n(T)\sim\tfrac{\sqrt{2n}}\pi$, so
$\Cal F_{\Pi}\ne\emptyset$ and $T$ is LLB, whence quasi finite.
\bul By Proposition 6.1,  if $\xi\subset\B$ is quasi finite,  then
$$\tfrac1{a_n(T)}\log\tfrac1{m(\xi_1^n(T)(x))}\overset{\goth
d}\to\rightarrow \u h(T,\xi)|\Cal N|\tag{\cancer}$$ as $n\to\infty$.
\heading{\tt\S7 Analogous properties of $\ppt$s }\endheading
\

The last part of this paper is devoted to the construction of an \f$\eppt$ having a generating partition with properties  analogous to (\cancer).
The "measure theoretic invariant" related to this is {\tt entropy dimension} as in [FP].
\sms Let $(\Bbb T,\Cal
T,m_{\Bbb T},R)$ be an irrational rotation of the circle (equipped
with Borel sets and Lebesgue measure).\sms Let $f\in L^2(\Bbb T)$
satisfy the weak invariance principle i.e. $B_n(t)\lra B(t)$  in
distribution on $C([0,1])$ where $B$ is Brownian motion and
$$B_n(t):=f_{[nt]-1}+(nt-[nt])f\circ T^{[nt]}$$ (where $f_k:=\sum_{j=0}^{k-1}f\circ R^j$).
Existence of such $f\in L^2(\Bbb T)$ is shown in [V].\bul In
particular, $$\tfrac{L_n}{\sqrt{n}},\ \tfrac{R_n}{\sqrt{n}}\ \ \
\overset{\goth d}\to\longrightarrow\ \ \   |\Cal N|,\ \ \tfrac{L_n+R_n}{\sqrt{n}}\ \ \
\overset{\goth d}\to\longrightarrow\ \ \   \Cal R$$ where
$R_n:=\max_{1\le k\le n}f_k$, $L_n:=\max_{1\le k\le n}(-f_k)$ and $\Cal R:=\max_{t\in [0,1]}B(t)-\min_{t\in [0,1]}B(t)$.
\

The  random variable $\Cal R$  is known as the {\it range of Brownian motion}. Its (non-Gaussian) distribution of is calculated in [Fe].
\sms Let $(Y,\C,\mu,S)$ be the $2$-shift with generating partition
$Q=\{Q_0,Q_1\}$ and symmetric product measure.

Let  $\rho:Y\to\Bbb R$ be defined by $\rho=\a_01_{Q_0}+\a_11_{Q_1}$ where
$\a_0<\a_1,\ \int_Y\rho d\mu=1$ and $\a_0,\ \a_1$ are rationally independent, then the {\it special
flow} (under $\rho$) $(Y^\rho,\C^\rho,q,S^\rho)$ is Bernoulli where
$$Y^\rho:=\{(y,s):\ y\in Y,\ s\in [0,\rho(y))\},\ \C^\rho:=\C\x\text{\rm Lebesgue},\ q:=\mu\x\l,$$ and
$$S^\rho_t(y,s):=(S^ny,s+t-\rho_n(y))$$ where $0\le s+t-\rho_n(y)<\rho(S^ny),\
\rho_n:=\sum_{j=0}^{n-1}\rho\circ S^j$. \bul Note that the ``{\tt
vertical}" partition $\o Q:=\{\o Q_0,\o Q_1\}$ where $\o Q_i:=Q_i\x
[0,\a_i)\ \ (i=0,1)$ generates $\C$ under $S^\rho$.

Define the \f$\ppt$ $(X,\B,m,T)$ by
$$X:=\Bbb T\x Y^\rho,\ m=m_{\Bbb T}\x q,\ \B:=\Cal T\x\C^\rho,\
T(x,(y,s)):=(R(x),S^\rho_{f(x)}(y,s)).\tag{\phone}$$

\

For $P$ a finite partition of $\Bbb T$ into intervals
(which generates $\Cal T$ under $R$), define the partition
$\xi=\xi_P$ of $X$ by
$$\xi(\om,y,s):=P(\om)\x \(\bigvee_{t\in\ \iota(0,f(\om))}S^\rho_{-t}\o Q\)(y,s)\tag{\uranus}$$
where for $x,y\in\ \Bbb R,\ \iota(x,y):=[x\wedge y,x\vee y]$ (the closed interval joining $x$ and $y$).
Next, we show that that $\xi$ is measurable and $H(\xi)<\infty$.
\proclaim{Proposition 7.1} \sms The partition \ \ $\xi$ is measurable,  generates
$\B$ under $T$, $H(\xi)<\infty$ and
$$\tfrac1{\sqrt n}I(\xi_0^{n-1}(T))\ \underset{n\to\infty}\to{\overset{\goth
d}\to{\lra}}\ h(S^\rho)\Cal R\tag{\neptune}$$ where $\Cal R$ is the range of Brownian motion.\endproclaim\sms{\bf Proof} The proof is in stages.
We claim first that
$$\xi_0^{n-1}(T)(\om,y,s)=P_0^{n-1}(R)(\om)\x\(\bigvee_{t\in [-L_n(\om),R_n(\om)]}S^\rho_{-t}t\o Q\)(y,s).\tag{\eighthnote}$$
\pf{\it\,of} (\eighthnote):   Note that for $n\ge 1$,
$$\align (T^{-n}\xi)&(\om,y,s) =\xi(R^n(\om),S_{f_n(\om)}^{\rho}(y,s))\\ &=P(R^n(\om))\x \(\bigvee_{t\in \iota(0,f(R^n(\om))}S^\rho_{-t}t\o Q\)(S_{f_n(\om)}^{\rho}(y,s))\\ &=
P(R^n(\om))\x \(\bigvee_{t\in \iota(f_n(\om),f_n(\om)+f(R^n(\om))}S^\rho_{-t}t\o Q\)(y,s)\\ &=
P(R^n(\om))\x \(\bigvee_{t\in \iota(f_n(\om),f_{n+1}(\om))}S^\rho_{-t}t\o Q\)(y,s)
.\endalign$$
To continue, we need the following (elementary) proposition:
\Par \ Let $a_n\in\Bbb R\ \ (n\ge 1)$ then
$\bigcup_{k=0}^{n-1}\iota(s_k,s_{k+1})=[m_n,M_n]$ where $a_0:=0$,
\f$s_n:=\sum_{k=0}^na_k,\ m_n:=\min_{0\le k\le n}s_k,\ M_n:=\max_{0\le k\le n}s_k.$

\

To finish the proof of (\eighthnote):
$$\align \xi_0^{n-1}(T)(\om,y,s)& =\bigvee_{k=0}^{n-1}T^{-k}\xi(\om,y,s)\\ &=
\bigcap_{k=0}^{n-1}P(R^k(\om))\x \(\bigvee_{t\in \iota(f_k(\om),f_{k+1}(\om))}S^\rho_{-t}t\o Q\)(y,s)\\ &=
P_0^{n-1}(R)(\om))\x \(\bigvee_{t\in \bigcup_{k=0}^{n-1}\iota(f_k(\om),f_{k+1}(\om))}S^\rho_{-t}t\o Q\)(y,s)\\ &\overset{\P}\to{=}
P_0^{n-1}(R)(\om)\x\(\bigvee_{t\in [-L_n(\om),R_n(\om)]}S^\rho_{-t}\o Q\)(y,s)
.\ \qed\ \text{(\eighthnote)}.\endalign$$

\par Now consider $\rho_n:Y\to\Bbb R$ defined by
$$\rho_n(y):=\cases & \sum_{k=0}^{n-1}\rho(S^ky)\ \ \ \ \ n>0,\\ & 0\ \ \ \ \ n>0,\\ & \sum_{k=1}^{|n|}\rho(S^{-k}y)\ \ \ \ \
n<0,\endcases$$ then $\rho_n(y)<\rho_{n+1}(y)$ and $\forall\ y\in Y$,
$\rho_n(y)\to\pm\infty$ as $n\to\pm\infty$.  \par For $y\in Y$ and
$t\in\Bbb R$, define $[t]_y\in\Bbb Z$ be so that $\rho_{[t]_y}(y)\le t<\rho_{[t]_y+1}(y)$.

It follows that for $t\in\Bbb R$:
\bul $\tfrac{|t|}{\a_1}-1\le |[t]_y|\le \tfrac{|t|}{\a_0}$; and
\bul $S^\rho_t(y,s)=(S^{[s+t]_y}y,s+t-\rho_{[s+t]_y}(y))$.

\

Our next claim is that
$$\xi_0^{n-1}(T)(\om,y,s)=P_0^{n-1}(R)(\om)\x Q_{[s-L_n(\om)]_y}^{[s+R_n(\om)]_y}(S)(y)\x \eta_n(\om,y)(s)\tag{\twonotes}$$
where for each $(\om,y)\in \Om\x Y,\ \eta_n(\om,y)$ is a partition of $[0,\rho(y))$ into at most $\tfrac{R_n(\om)+L_n(\om)+1}{\a_0}$ intervals.

\pf{\it of} (\twonotes).

Fixing $(\om,y,s)\in X$ and $n\ge 1$, we have

$$\align \(\bigvee_{t\in [-L_n(\om),R_n(\om)]}S^\rho_{-t}\o Q\)(y,s)&=
\bigcap_{t\in [-L_n(\om),R_n(\om)]}\o Q(S^\rho_{t}(y,s))\\ &=
\bigcap_{t\in [-L_n(\om),R_n(\om)]}Q(S^{[s+t]_y}y)\x [0,\rho(S^{[s+t]_y}y))\\ &=
\bigcap_{j\in [[s-L_n(\om)]_y,[s+R_n(\om)]_y]}S^{-j}Q(y)\x \eta_n(\om,y,s)\\ &=
Q_{[s-L_n(\om)]_y}^{[s+R_n(\om)]_y}(S)(y)\x \eta_n(\om,y)(s).\endalign$$
where for each $(\om,y)\in \Om\x Y,\ \eta_n(\om,y)$ is a partition of $[0,\rho(y))$ into at most $[R_n(\om)]_y]-[-L_n(\om)]_y\le \tfrac{R_n(\om)+L_n(\om)+1}{\a_0}$ intervals.\ \qed\ (\twonotes).

\

\bul Observation of (\twonotes) with $n=1$ shows that
$$\xi(\om,y,s):=P(\om)\x Q_{-\nu_-(\om,y,s)}^{\nu_+(\om,y,s)}(S)(y)\x
\eta_1(\om,y)(s)$$ where
$$\nu_{+}(w,y,s)=[s+f(\om)\vee 0]_y,\ \
\nu_{-}(w,y,s)=[s+f(\om)\wedge 0]_y.$$ Thus, $\xi$ is measurable.

\sms Moreover, writing $\Cal Z:=\{[\nu_-=k,\ \nu_+=\ell]:\ k,\ell\in\Bbb Z\}$, we see that

$$\align I(\xi|\Cal Z)(\om,y,s)&=I(P)(\om)+I(Q_{[s+f(\om)\wedge 0]_y}^{[s+f(\om)\vee 0]_y})(S)(y)+I(\eta_1(\om,y)(s)\\ &\le I(P)(\om)+([s+f(\om)\wedge 0]_y+[s+f(\om)\vee 0]_y)\cdot\log 2+\log\tfrac{1+|f(\om)|}{\a_0}\\ &\le I(P)(\om)+\tfrac{|f(\om)|+1}{\a_0}\cdot\log 2
+\log\tfrac{1+|f(\om)|}{\a_0}\endalign$$
and
$$H(\xi|\Cal Z)\le H(P)+\tfrac{\log 2}{\a_0}(\|f\|_1+1)+\int_{\Om}\log\tfrac{1+|f|}{\a_0}dm<\infty.$$

Now $|\nu_\pm(\om,y,s)|\le \tfrac{|f(\om)|+1}{\a_0}$ and
$$(\nu_{+}(w,y,s),\nu_{-}(w,y,s))=\cases & ([s+f(\om)\vee 0]_y,0)\ \ \ \ \ f(\om)\ge 0,\\ &
(0,[s+f(\om)\wedge 0]_y)\ \ \ \ \ f(\om)< 0;\endcases$$
whence using ($\bigstar$) (see page 6) $H(\Cal Z) <\infty$ and
$$H(\xi)=H(\xi|\Cal Z)+H(\Cal Z) <\infty.$$

\

\bul Since $\xi$ is measurable, (\eighthnote) now shows that it generates
$\B$ under $T$.
\bul To establish (\neptune), we claim that for a.e.
$(x,y,s)$, for
 any $\e>0$, for sufficiently large
$n=n(x,y,s)$,
$$\align P_0^{n-1}(R) &(x)\x Q_{-L_n(x)(1+\e)}^{R_n(x)(1+\e)}(S)(y)\x\eta_n(x,y)(s)\subseteq\xi_0^{n-1}(T)(x,y,s)\\
 &\subseteq P_0^{n-1}(R)(x)\x Q_{-L_n(x)(1-\e)}^{R_n(x)(1-\e)}(S)(y)\x \eta_n(x,y)(s)\tag{$\clubsuit$}\endalign$$
where for each $(\om,y)\in \Om\x Y,\ \eta_n(\om,y)$ is a partition of $[0,\rho(y))$ into at most $\tfrac{R_n(\om)+L_n(\om)+1}{\a_0}$ intervals.

\pf {of} ($\clubsuit$): For a.e. $(x,y,s)\in X$, $R_n(x),\ L_n(x)\uparrow\ \infty$ and $\rho_n(y)\sim n$, whence $|[s-L_n(x)]_y|\sim L_n(\om)$ and $[s+R_n(\om)]_y\sim R_n(x)$.
($\clubsuit$) follows from (\twonotes) using this.\qed

\bul We claim next
that $\forall\ (x,y)\in\Bbb T\x Y$,
$$\tfrac1{\sqrt n}(I(P_0^{n-1}(R))+
I(\eta_n(x,y))\overset{m}\to\longrightarrow 0.\tag{\sun}$$
 \pf $\#\eta_n(x,y)\le \Cal E_n(x):=\tfrac{R_n(x)+L_n(x)+1}{\a_0}$ and $\#P_0^{n-1}(R)\le Mn$ for some $M>0,\ \forall\ n\ge 1$, whence

$$m([I(P_0^{n-1}(R))\ge t\sqrt n])
 \le \tfrac1{t\sqrt n}H(P_0^{n-1}(R))\lesssim\tfrac{\log n}{t\sqrt n}\ \to 0\ \text{\tt as}\
n\to\infty$$

and $\forall\ (x,y)$,

$$m([I(\eta_n(x,y)(s))])\ge t\sqrt n])
\le \tfrac1{t\sqrt n}H(\eta_n(x,y))\le\tfrac{\log\Cal E_n(x)}{t\sqrt n}\overset{m}\to\lra 0\ \text{\tt as}\
n\to\infty$$ proving (\sun).\ \checkmark

\

Using ($\clubsuit$), (\sun) and ($\goth I$) for $S$
 we have, as $n\to\infty$, $$\align \tfrac1{\sqrt n}I(\xi_0^{n-1}(T))(x,y,s) &= \tfrac1{\sqrt n}I({Q_{ -L_n(x)(1+o(1))}^{R_n(x)(1+o(1))}(S)})(y)+O(\tfrac{\log
n}{\sqrt n})\\ &= \tfrac1{\sqrt n}(L_n(x)+R_n(x))\log 2(1+o(1))+O(\tfrac{\log
n}{\sqrt n})\\ &\overset{\goth
d}\to{\lra}\Cal R \log 2\\ &=\Cal R h(S^\rho ).\qed\ (\neptune)\endalign$$
\subheading{Estimation of entropy dimension}

\

Let $(Z,\Cal D,\nu,R)$ be a $\ppt$ and let $P\subset \Cal D$ be a countable partition of $Z$.

As in  [FP], let for $n\ge 1,\ \e>0,\ a=\bigcap_{k=0}^{n-1}R^{-k}a_k\in P_0^{n-1}(R)$,
$$B(n,P,a,\e):=\bigcup_{a'\in P_0^{n-1}(R),\ \o d(a,a')<\e}a$$
where $\o d(a,a'):=\tfrac1n\#\{0\le k\le n-1:\ a_k\ne a_k'\}$ is Hamming distance, and let
$$K(P,n,\e):=\min\,\{\# F:\ F\subset P_0^{n-1}(R),\ \nu(\bigcup_{a\in F}B(n,P,a,\e))>1-\e\}.$$
\bul The $\eppt$ is said to have {\it upper entropy dimension} $\Delta\in [0,1]$ if for some countable, measurable generating partition
$P$ with finite entropy (and hence  -- as proved in [FP]-- for all such),
$$ \varlimsup_{n\to\infty}\tfrac{\log\log K(P,n,\e)}{\log n} \underset{\e\to 0}\to\longrightarrow\ \Delta.$$

\proclaim{Proposition 7.2} Let $(X,\B,m,T)$ be as in (\phone), then the upper entropy dimension is at most $\frac12$.\endproclaim\pf
\ Let $\xi=\xi_P$ be as in (\uranus) and let $h=h(S^\rho)$. For $n\ge 1,\ J\subset \Bbb R_+$ an  interval bounded away from $0$ and $\infty$,  define
$\xi_n(J):=\{a\in\xi_0^{n-1}(T):\ \tfrac1{\sqrt n}\log\tfrac1{m(a)}\in hJ\}$.
\Par \ \   We claim that
$\#\xi_n(J)\sim E(1_J(\Cal R)e^{h\Cal R\sqrt n})e^{o(\sqrt n)}$ as $n\to\infty$.
\pf Suppose that $J=[r-\d,r+\d]$, then
$$\align P(\Cal R\in J)\ & \lla\ m([\tfrac1{\sqrt n}I(\xi_0^{n-1}(T))\in hJ])\\ &=\sum_{a\in\xi_n(J)}m(a)\\ &=\#\xi_n(J)e^{-h\sqrt n (r\pm\d)}\endalign$$
(because $m(a)=e^{-h\sqrt n (r\pm\d)}\ \forall\ a\in\xi_n(J)$);\ whence
$$E(e^{h\sqrt n (\Cal R-2\d)}1_J(\Cal R))\lesssim \#\xi_n(J)\lesssim E(e^{h\sqrt n (\Cal R+2\d)}1_J(\Cal R)).$$
Using this on a decomposition of  $J$ into a finite union of disjoint short enough intervals  proves $\#\xi_n(J)=E(e^{h\sqrt n \Cal R}1_J(\Cal R))e^{\pm\e\sqrt n}\ \forall\ \e>0$, whence \P.\qed
\bul Evidently $K(\xi,n,\e)\le\#\xi_n([\tfrac1M,M])$ for some $M=M_\e>0$ whence
$K(\xi,n,\e)\le e^{c_\e\sqrt n(1+o(1))}$ and $\varlimsup_{n\to\infty}\tfrac{\log\log K(\xi,n,\e)}{\log n}\le\tfrac12\ \forall\ \e>0$.\ \ \Box
\subheading{\tt Remark on the lower bound} 
\

The upper estimate for the entropy dimension  follows from the the weak invariance principle for the ``{\tt random walk}" $f_n$. In a similar manner, a lower estimate would follow from an analogous result for the ``{\tt local time}" of the random walk. Such a result is not available for the present example. However, such considerations show that the ``{\tt relative entropy dimension}" of an aperiodic, centered random walk in random scenery  over its Bernoulli factor is $1/2$. 

 \heading References\endheading 
 \enddocument